\documentclass[11pt]{article}
\usepackage{latexsym,amsthm,amsmath,amssymb,graphicx}
\usepackage{todonotes}
\usepackage[utf8x]{inputenc}
\usepackage{hyperref}
\usepackage{soul}
\usepackage{authblk}

\newtheorem{theorem}{Theorem}
\newtheorem{lemma}{Lemma}

\newtheorem{corollary}{Corollary}

\newtheorem{remark}{Remark}

\newcommand{\codes}{{\textsc{Codes}}}

\newcommand{\x}{\gamma^{\rm{X}}}
\newcommand{\id}{\gamma^{\rm{ID}}}
\newcommand{\itd}{\gamma^{\rm{ITD}}}
\newcommand{\fd}{\gamma^{\rm{FD}}}
\newcommand{\ftd}{\gamma^{\rm{FTD}}}
\newcommand{\ld}{\gamma^{\rm{LD}}}
\newcommand{\ltd}{\gamma^{\rm{LTD}}}
\newcommand{\od}{\gamma^{\rm{OD}}}
\newcommand{\otd}{\gamma^{\rm{OTD}}}

\newcommand{\calN}{\mathcal{N}}

\newcommand{\dc}[1]{{\color{teal} #1}}

\makeatletter
\if@todonotes@disabled

\else

\fi
\makeatother

\title{On lower bounds for cardinalities of several separating-dominating codes in graphs}

\author[1,2]{Dipayan Chakraborty}
\author[1]{Annegret K. Wagler}

\affil[1]{LIMOS, Universit\'{e} Clermont Auvergne, CNRS, Mines Saint-\'{E}tienne, Clermont Auvergne INP, LIMOS, 63000 Clermont-Ferrand, France.}
\affil[2]{Department of Mathematics and Applied Mathematics, University of Johannesburg, Auckland Park, 2006, South Africa}

\providecommand{\keywords}[1]
{
  \small	
  \textbf{\textit{Keywords---}} #1
}

\begin{document}

\maketitle

\begin{abstract}
 In the literature, several different identification problems in graphs have been studied, the most widely studied such problems are the ones based on dominating sets as a tool of identification. Hereby, the objective is to separate any two vertices of a graph by their unique neighborhoods in a suitably chosen dominating or total-dominating set. Such a (total-)dominating set endowed with a separation property is often referred to as a \emph{code} of the graph. The problems of determining such codes of minimum cardinality are all shown to be NP-hard. A typical line to attack such problems is, therefore, to provide bounds on the cardinalities of the studied codes.
  
  In this paper, we are interested in extremal graphs for lower bounds in terms of the order of the graph.
  For some codes, logarithmic lower bounds are known from the literature. 
  We provide for eight different identification problems a general construction of extremal graphs minimizing the cardinality of a minimum code compared to the order of the graph. This enables us to
  reprove the existing logarithmic lower bounds,
  to establish such bounds for further codes, and to characterize all graphs attaining these bounds.
\end{abstract}


\section{Introduction}

\begin{figure}[!t]
\begin{center}
\includegraphics[scale=1.0]{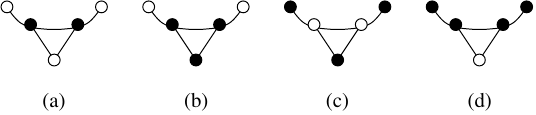}
\caption{Minimum X-codes in a graph (the black vertices belong to the code), where (a) is both an LD- and LTD-code, (b) both an OD- and OTD-code, (c) an ID-code, (d) an ITD-code as well as both an FD- and FTD-code.}
\label{fig_exp_X-codes}
\end{center}
\end{figure}

In the area of identification problems in graphs, the aim is to distinguish any two vertices of a graph by the unique intersection of their neighborhoods with a suitably chosen dominating or total-dominating set of the graph. Consider a simple, undirected graph $G=(V,E)$ and denote by $N(v)$ the \emph{open neighborhood} of a vertex $v \in V$. Moreover, let $N[v]$ denote the \emph{closed neighborhood} of $v$, that is, $N[v] = N(v) \cup \{v\}$. A subset $C \subseteq V$ is \emph{dominating} (resp. \emph{total-dominating}) if $N[v]\cap C$ (resp. $N(v)\cap C$) are non-empty sets for all $v \in V$. 
To distinguish vertices of a graph by means of a (total-)dominating set, different \emph{separation} properties of these sets have been studied.
A subset $C \subseteq V$ is
\begin{itemize}
  \itemsep-3pt
\item \emph{location-separating}\footnote{Location-separating sets are typically called \emph{locating sets} in the literature, see e.g. \cite{HHH_2006,S_1988}. For consistency in terminology, we call them here location-separating sets to specify that they refer to one of the studied separation properties.} (for short \emph{L-separating}) if $N(v)\cap C$ is a unique set for each $v \in V\setminus C$, 
\item \emph{open-separating} (for short \emph{O-separating}) if $N(v)\cap C$ is a unique set for each $v \in V$,
\item \emph{closed-separating} (for short \emph{I-separating}) if $N[v]\cap C$ is a unique set for each $v \in V$, and
\item \emph{full-separating} (for short \emph{F-separating}) if it is both open- and closed-separating, that is if $N[u] \cap C \ne N[v] \cap C$ and $N(u) \cap C \ne N(v) \cap C$ holds for each pair of distinct vertices $u,v \in V$.
\end{itemize}
Combining one separation property with one domination property, the following identification problems have been studied in the literature:
\begin{itemize}
  \itemsep-3pt
\item location-separation with domination and total-domination leading to \emph{locating dominating codes} (for short \emph{LD-codes}) and \emph{locating total-dominating codes} (for short \emph{LTD-codes}), see \cite{S_1988} and \cite{HHH_2006}, respectively;
\item open-separation with domination and total-domination leading to \emph{open-sepa\-ra\-ting dominating codes} (\emph{OD-codes} for short) and \emph{open-separating total-do\-minating codes}\footnote{OTD-codes were introduced independently in~\cite{HLR_2002} and in~\cite{SS_2010} under the names of \emph{strongly $(t,\le l)$-identifying codes} and \emph{open neighborhood locating-dominating sets} (or \emph{OLD-sets}), respectively. However, due to consistency in naming that specifies the separation and the domination property, we prefer to call them open-separating total-dominating codes in this article.
} (\emph{OTD-codes} for short), see \cite{CW_ISCO2024} and  \cite{SS_2010}, respectively;
\item closed-separation with domination and total-domination leading to \emph{identifying codes} (\emph{ID-codes} for short) and \emph{identifying total-dominating codes}\footnote{Identifying total-dominating codes had been introduced to the literature in~\cite{HHH_2006} under the name \emph{differentiating total-dominating codes}. However, due to consistency in notation, we prefer to call them ITD-codes in this article.} (\emph{ITD-codes} for short), see \cite{KCL_1998} and \cite{HHH_2006}, respectively;
\item full-separation with domination and total-domination leading to \emph{full-sepa\-ra\-ting dominating codes} (for short \emph{FD-codes}) and \emph{full-separating total-do\-minating codes} (for short \emph{FTD-codes}), respectively, see \cite{CW_2024-FS}.
\end{itemize}
Figure \ref{fig_exp_X-codes} illustrates examples of such codes in a small graph. 
%%%%%%%%%%%%%%
%%%%%%%%%%%%%%

Codes of these types have manifold applications, for example, in locating intruders in facilities using sensor networks~\cite{UTS_2004}. To model such monitoring problems, 
the rooms of a building are represented by vertices of a graph, doors between rooms as its edges, and the task is to place sensors at some vertices so that the intruder (like a fire, a thief, or a saboteur) can be either directly detected in the room equipped with a sensor or can be located with the help of a unique pattern of alerts sent out by sensors installed in some of the adjacent rooms.
Hereby, the subset of rooms equipped with sensors corresponds to a code in the graph. The domination properties guarantee that all rooms are monitored, whereas the separation properties ensure that the intruder can be exactly located. 
For the above described setting, location-separation is used \cite{S_1988}.
Maintaining the functionality if some sensors cannot distinguish between the presence of an intruder in the room itself or in an adjacent room corresponds to closed-separation \cite{KCL_1998},
monitoring with partly disabled sensors leads to open-separation \cite{HLR_2002,SS_2010}, 
whereas full-separation can handle both fault types simultaneously \cite{CW_2024-FS}.

It is known from the literature that not every graph has an X-code, where $X \in \codes$ = \{LD, LTD, OD, OTD, ID, ITD, FD, FTD\}.
Summarizing results from for instance \cite{CW_ISCO2024,CW_2024-FS,HHH_2006,KCL_1998,SS_2010}, we see that a graph $G$ has no X-code involving
%%%%%%%%%%%%%%%%%%%%%%
\begin{itemize}
  \itemsep-3pt
\item total-domination if $G$ has \emph{isolated vertices}, that is, vertices $v$ with $N(v) = \emptyset$;
\item open-separation if $G$ has \emph{open twins}, that is, non-adjacent vertices $u,v$ with $N(u) = N(v)$;
\item closed-separation if $G$ has \emph{closed twins}, that is, adjacent vertices $u,v$ with $N[u] = N[v]$.
\end{itemize}

A graph without open twins (resp. closed twins) is called \emph{open-twin-free} (resp. \emph{closed-twin-free}). Moreover, a graph that is both open-twin-free and closed-twin-free is called \emph{twin-free}. In addition, a graph without isolated vertices is called \emph{isolate-free}. Table \ref{tab_admissible} illustrates which X-codes are prevented from existing in a graph as a result of the presence of isolated vertices, open twins and closed twins in the graph.
\begin{table}[ht]
\begin{center}
\begin{tabular}{ r || c | c || c | c || c | c || c | c }
 X                 & LD   & LTD  & OD   & OTD  & ID & ITD & FD & FTD \\ \hline
 isolated vertices &      &  x   &      &  x   &      &  x   &      &  x   \\
 open twins        &      &      &  x   &  x   &      &      &  x   &  x   \\
 closed twins      &      &      &      &      &  x   &  x   &  x   &  x  \\
\end{tabular}
\end{center}
\caption{
The X-codes that do not exist in graphs having isolated vertices or open or closed twins are marked with an x in the table.}
\label{tab_admissible}
\end{table}
%%%%%%%%%%%%%%%%%%%%%%

Calling a graph $G$ \emph{X-admis\-sible} if $G$ has an X-code, we see, for example, from Table \ref{tab_admissible} that every graph $G$ is LD-admissible, whereas a graph $G$ must neither have isolated vertices nor open or closed twins to be FTD-admissible.

Given an X-admis\-sible graph $G$, the X-problem on $G$ is the problem of finding an X-code of minimum cardinality $\x(G)$ in $G$, called its \emph{X-number}. 
Problems of this type have been actively studied in the context of various applications during the last decade, see e.g. the bibliography maintained by Jean and Lobstein \cite{Lobstein_Lib}. 
In particular, all X-problems for any $X \in \codes$ have been shown to be NP-hard \cite{CW_ISCO2024,CW_2024-FS,CHL_2003,CSS_1987,SS_2010}.
Therefore, a typical line to attack such problems is to determine closed formulas for special graphs or to provide bounds on the X-numbers.
We again refer to \cite{Lobstein_Lib} for a literature overview. 

%%%%%%%%%%%%%%%%%%%%%%%%
\begin{figure}[h]
\begin{center}
\includegraphics[scale=0.8]{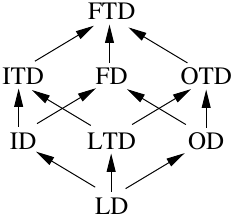}
\caption{The relations between the X-numbers for all $X \in \codes$, where $X' \longrightarrow X$ stands for $\gamma^{X'}(G) \leq \gamma^X(G)$.}
\label{fig_relations}
\end{center}
\end{figure}
%%%%%%%%%%%%%%%%%%%%%%%%
A general relation between the X-numbers for all $X \in \codes$ has been established in \cite{CW_2024-FS}, see Figure \ref{fig_relations}.
We see in particular that $\ld(G)$ is a lower bound for \emph{all} other X-numbers, whereas $\ftd(G)$ is an upper bound for \emph{all} other X-numbers in FTD-admissible graphs. Moreover, it has been shown in \cite{CW_ISCO2024} (resp. \cite{CW_2024-FS}) that the OD- and OTD-numbers (resp. the FD- and FTD-numbers) of any OTD-admissible (resp. FTD-admissible) graph differ by at most one. 

In this paper, we are interested in extremal graphs for lower bounds on X-numbers in terms of the order of the graph.
As the separation properties require distinct intersections of the neighborhoods with an X-code, it is immediate that such bounds are logarithmic in the order of the graph.
For X-admissible graphs $G$ of order $n$, it is known so far that the following bounds hold:
\begin{flalign} \label{eq:LB_literature}
\x(G) \geq \left\{
\begin{array}{ll}
  \lfloor \log n \rfloor   & \mbox{if } X \in \{LD,LTD\} \ \cite{HR_2012,S_1988} \\
  \lceil \log n \rceil                   & \mbox{if } X = OD \ \cite{CW_ISCO2024} \\
  \lceil \log (n+1) \rceil & \mbox{if } X \in \{ID,OTD\} \ \cite{KCL_1998,SS_2010} \\
\end{array}
\right.
\end{flalign}

All logarithms in the article are to the base~$2$. Examples of graphs whose X-numbers attain these logarithmic lower bounds have already been discussed in several articles in the literature. For LD-codes, a characterization of graphs whose LD-numbers attain the logarithmic lower bound in (\ref{eq:LB_literature}) is given in~\cite{S_1988}. For LTD-codes studied in~\cite{HR_2012}, the authors lay out several examples (but without characterization) of graphs $G$ with their LTD-numbers being equal to the lower bound in (\ref{eq:LB_literature}). For ID-codes, again a complete characterization of all graphs with their ID-numbers equal to the lower bound in (\ref{eq:LB_literature}) was presented in~\cite{M_2006}. In~\cite{SS_2010}, one of the papers introducing the OTD-codes, the authors discuss certain examples (but again without characterization) of graphs whose OTD-numbers attain the respective logarithmic lower bound in (\ref{eq:LB_literature}). 
In~\cite{CW_ISCO2024} as well, the lower bound of OD-numbers in (\ref{eq:LB_literature}) was introduced without characterization of the graphs realizing such bounds.

To the best of our knowledge, the above works comprise all that have been studied in the literature of identification problems in graphs with respect to general lower bounds on the X-numbers of graphs and the characterization of such graphs realizing these bounds. As can be seen, while a characterization of these extremal graphs exist only for LD- and ID-codes~\cite{S_1988, M_2006}, such characterizations are lacking for the other codes that have already long existed in the literature. Moreover, in the recent past, three new codes, namely OD-, FD- and FTD-codes, have been introduced to the literature (see~\cite{CW_ISCO2024, CW_2024-FS}) which too so far lack a study dedicated to the extremal graphs attaining such lower bounds. Considering this state of the art, our aim is to study, for all $X \in \codes$, the characterization of all extremal graphs whose X-numbers attain the corresponding logarithmic lower bound in terms of their order. To that end, we provide, for all $X \in \codes$, a general construction technique for graphs minimizing their X-number compared to their order. This allows us to reprove the results stated in (\ref{eq:LB_literature}) (this time under the common construction technique) and to extend
the list in (\ref{eq:LB_literature}) with the
logarithmic lower bounds for the remaining X-numbers with $X \in \{ITD,FD,FTD\}$. Moreover, we characterize for all $X \in \codes$ the graphs attaining these logarithmic lower bounds. In addition, our work implies a small correction to a result in~\cite{S_1988} regarding the number of maximum-ordered graphs with a common LD-number. We close with some concluding remarks and lines of future research.

\section{Extremal graphs for lower bounds on X-numbers}
\label{Sec_2}

The aim of this section is to construct extremal graphs minimizing their X-number with respect to their order. To that end, we define the \emph{X-density} of an X-admissible graph $G$ on $n$ vertices to be the ratio $\frac{\x(G)}{n}$. Hence, in this section, we look for graphs with the minimum possible X-densities.

Given any set $C$, the set of all subsets of $C$ is called the \emph{power set} of $C$ and is denoted by $2^C$. Then, $|2^C| = 2^{|C|}$. As the separation properties require unique intersections of the neighborhoods of vertices with an X-code $C$,
we call $N(v)\cap C$ (resp. $N[v]\cap C$) the \emph{open} (resp. \emph{closed}) \emph{signature} of a vertex $v$. Therefore, notice that the signature of a vertex is one out of the $2^{|C|}$ possible subsets of $C$. Moreover, given a code $C$ of a graph, the signatures of the vertices of $C$ themselves, in particular, play an important role in several counting arguments later. As a result, for a vertex subset $C$ of a graph $G$, we define $\calN(C) = \{N(v) \cap C : v \in C\}$ and $\calN[C] = \{N[v] \cap C : v \in C\}$. In addition, let $\calN_C = \calN(C) \cup \calN[C]$.

\begin{remark} \label{rem:G[C] X-admissible}
Let $X \in \codes$ and let $G$ be an X-admissible graph. Moreover, let $C$ be an X-code of $G$. Then the subgraph $G[C]$ is X-admissible.
\end{remark}

In the following, we provide a common construction for extremal graphs $G^S(k) = (C \cup V \setminus C, E)$ for the four studied separation properties $S \in \{L,O,I,F\}$ such that $C$ is a minimum X-code of cardinality $k$ and the order $|V|$ of $G^S(k)$ is maximum possible. 
As no graph can have an X-code of cardinality 1 for $X \in \{LTD, OTD, ITD, FTD\}$, we assume $k \geq 2$ in the sequel. In all four cases $S \in \{L,O,I,F\}$, we require that $C$ is an $S$-separating set of $G^S[C]$ and that the maximum possible number of signatures is used, but we do not impose further conditions for the induced subgraphs $G^S(k)[C]$ and $G^S(k)[V \setminus C]$. Thus, even though a different edge set of $G^S(k)[C]$ or of $G^S(k)[V \setminus C]$ gives rise to a different graph $G^S(k)$, we simply denote all such graphs by $G^S(k)$ for brevity and choose the edge sets of $G^S(k)[C]$ and $G^S(k)[V \setminus C]$ only as and when necessary.

In order to count the number of graphs $G^S(k)[C]$, we further introduce the following notations. For $S \in \{L,O,I,F\}$, let $\eta^S(k)$ (respectively $\overline{\eta}^S(k)$) denote the number of all (respectively isolate-free) graphs on $k$ vertices that admit an $S$-separating set. Then the number of graphs on $k$ vertices with an isolated vertex and admitting an $S$-separating set is $\overline{\eta}^S(k-1)$. In addition, let $\eta(k)$ denote the total number of graphs on $k$ vertices. In other words, we have $\eta(k) = 2^{k \choose 2}$. In particular, since any graph admits a locating set (refer to Table~\ref{tab_admissible}), we therefore have $\eta^L(k) = \eta(k) = 2^{k \choose 2}$.

\subsection{The case of LD- and LTD-codes}
\label{Sub_2_L}

In this subsection, we construct graphs minimizing their X-density for $X \in \{LD, LTD\}$. Recall that for location-separating sets $C$, it is required that the open signatures $N(v)\cap C$ are unique sets for all vertices $v$ outside $C$. The following result is a combination of the works in~\cite{HR_2012} and~\cite{S_1988}. Nevertheless, we prove them here as a single result for the flow of arguments and the self-containment of this article.

\begin{lemma}\label{rem:L-cardinality}
For $X \in \{LD, LTD\}$, if $G$ is an X-admissible graph on $n$ vertices and with $\x(G) = k$, then $n \le 2^k-1+k$. In particular, therefore, we have $\x(G) \ge \lfloor \log n \rfloor$.
\end{lemma}

\begin{proof}
Let $C$ be a minimum X-code of $G$. Since at most $2^{|C|}-1$ non-empty subsets of $C$ can be open signatures of vertices in $V(G) \setminus C$, we have $|V(G) \setminus C| \le 2^{|C|}-1$. This implies $n = |V(G)| = |C| + |V \setminus C| \le 2^{|C|}-1+|C| = 2^k-1+k$.

To prove the second result, we observe that $2^k-1+k < 2^{k+1}$ for all $k \ge 1$. This implies that $k > -1+\log n$ and hence, $k \ge \lfloor \log n \rfloor$.
\end{proof}

Next, we construct graphs realizing the minimum X-densities for $X \in \{LD, LTD\}$.
Define the graph $G^L(k) = (C \cup V \setminus C, E)$ for $k \geq 2$ as follows.
\begin{itemize}
\itemsep0pt
  \item $C$ is a vertex subset of $G$ of cardinality $k$,
  \item $V \setminus C$ has exactly one vertex $v_{C'}$ with signature $N(v_{C'}) \cap C = C'$ for each non-empty subset $C' \subseteq C$.
\end{itemize}

\begin{figure}
\begin{center}
\includegraphics[scale=1.0]{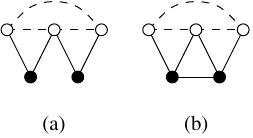}
\caption{The graphs $G^L(2)$. Black-filled vertices belong to the minimum code $C$, the dotted lines between un-filled vertices are optional. The graph in (a) is extremal for LD-codes, the graph in (b) is extremal for both LD- and LTD-codes.}
\label{fig_L}
\end{center}
\end{figure}

Notice that this construction imposes conditions only on those edges concerning the open signatures (with respect to $C$) of the vertices of $V \setminus C$, but not for the induced subgraphs $G^L(k)[C]$ and $G^L(k)[V \setminus C]$. As an example, Figure \ref{fig_L} depicts the possible graphs $G^L(k)$ for $k = 2$.

\begin{lemma} \label{rem:GL(k)-cardinality}
For the graph $G^L(k) = (C \cup V \setminus C, E)$, we have $|V| = 2^k-1+k$. Moreover, the number of graphs $G^L(k)$ that exist is the number of graphs on $k$ vertices times the number of graphs on $2^k-1$ vertices, that is $\eta(k) \times \eta(2^k-1)$.
\end{lemma}

\begin{proof}
By the construction of $G^L(k)$, each non-empty subset of $C$ is an open signature of a vertex in $V \setminus C$. Therefore, $|V \setminus C| = 2^{|C|}-1$ and hence, $|V| = |C|+|V \setminus C| = 2^{|C|}-1+|C| = 2^k-1+k$.

The subgraphs $G^L(k)[C]$ and $G^L(k)[V \setminus C]$ are on $k$ vertices and $2^k-1$ vertices, respectively. Therefore, the result on the number of graphs $G^L(k)$ follows by the fact that the construction of $G^L(k)$ has no restrictions on the edge sets of $G^L(k)[C]$ and $G^L(k)[V \setminus C]$.
\end{proof}

We next show that, for $X \in \{LD, LTD\}$, the graphs $G^L(k)$ are the ones with the minimum possible X-densities.

\begin{theorem}\label{thm_L}
Let $k \ge 2$. Then in $G^L(k) = (C \cup V \setminus C, E)$, $C$ is a minimum
\begin{itemize}
\itemsep0pt
  \item LD-code,
  \item LTD-code if $G^L(k)[C]$ has no isolated vertex, that is, if $G^L(k)[C]$ is LTD-admissible.
\end{itemize}
Moreover, for $X \in \{LD, LTD\}$, the minimum X-density of a graph with X-number $k$ is attained by any $G^L(k)$ for $X=LD$ and by those $G^L(k)$ with isolate-free $G^L(k)[C]$ for $X=LTD$.
\end{theorem}

\begin{proof}
The second condition in the construction of $G^L(k)$ ensures that each vertex in $V \setminus C$ has a unique open signature in $C$. Therefore, $C$ is a location-separating set of $G^L(k)$. Moreover, since no vertex in $V \setminus C$ has an empty open signature in $C$, the set $C$ total-dominates all vertices of $V \setminus C$. This implies that $C$ is a dominating set (and hence, an LD-code) of $G^L(k)$ and is a total-dominating set (and hence, and LTD-code) of $G^L(k)$ if $G^L(k)[C]$ has no isolated vertices. Moreover, in both the cases, the code $C$ must be of minimum cardinality, as otherwise, by Lemma~\ref{rem:L-cardinality}, we would have $|V| \le 2^{k-1}-1+(k-1) = 2^{k-1}-2+k$ contradicting $|V| = 2^k-1+k$ (by Lemma~\ref{rem:GL(k)-cardinality} and the fact that $2^{k-1}-2 < 2^k-1$ for all $k \geq 2$). 

Therefore, we have $\x(G^L(k)) = k$ both for $X=LD$ and for $X=LTD$ with isolate-free $G^L(k)[C]$. Thus, for $X \in \{LD, LTD\}$, the result that such graphs $G^L(k)$ have the minimum X-densities among all graphs with X-number $k$ follows by the comparison of graph orders in Lemmas~\ref{rem:L-cardinality} and~\ref{rem:GL(k)-cardinality}.
\end{proof}

Next, for $X \in \{LD, LTD\}$, we characterize all X-admissible graphs which attain the logarithmic lower bound in Theorem~\ref{rem:L-cardinality}.

\begin{theorem} \label{pro_L}
For $X \in \{LD, LTD\}$, let $G$ be an X-admissible graph of order $n$ and let $k = \lfloor \log n \rfloor \ge 2$. Then we have $\x(G) = k$ if and only if $G$ can be obtained from $G^L(k) = (C \cup V\setminus C, E)$ with X-admissible $G^L(k)[C]$ by removing up to $k-1$ vertices of $V \setminus C$.
\end{theorem}

\begin{proof}
To prove the sufficiency part, let $G \cong G^L(k) - S$ for some graph $G^L(k)$, where $G^L(k)[C]$ is X-admissible and $S \subset V \setminus C$ with $|S| \le k-1$. Then, by Theorem~\ref{thm_L}, the set $C$ is a minimum X-code of $G^L(k)$. Thus, $C$ is also an X-code of $G^L(k)-S$. Let $D \subset V(G)$ be the isomorphic copy of $C$ in $G$. Then, by the isomorphism, $D$ is an X-code of $G$. Moreover, we have $|D| = |C| = k = \lfloor \log n \rfloor$. Hence, by Lemma~\ref{rem:L-cardinality}, we have $\x(G) = k$.

To prove the necessary part, let us assume that $\x(G) = k$. Let $C$ be any minimum X-code of $G$. Then, it is clear that every vertex in $V(G) \setminus C$ has a unique open signature with respect to $C$. However, since the vertices in $V(G) \setminus C$ may not correspond to all the possible subsets of $C$ as open signatures, we therefore have  $G \cong G^L(k) - S$, where $S \subset V \setminus C$ and the edge set of $G^L(k)[V \setminus C]$ is so chosen to allow for the isomorphism. Also, $G^L(k)[C]$ is X-admissible since $C$ is an X-code of $G$. We now show that $|S| \le k-1$. Let us assume to the contrary that $|S| \ge k$. By Lemma~\ref{rem:GL(k)-cardinality}, we have $|V|= 2^k-1+k$. Therefore, we have $|V|-|S| \le 2^k-1 < 2^k$. Since $n = |V| - |S|$ by the isomorphism $G \cong G^L(k) - S$, we have $\lfloor \log n \rfloor = \lfloor \log (|V|-|S|) \rfloor < k$, a contradiction to $k = \lfloor \log n \rfloor$. This proves the necessary part.
\end{proof}

For instance, the graph $G^L(2)$ and all graphs obtained by removing any arbitrary un-filled vertex from $G^L(2)$ in Figure \ref{fig_L} (resp. Figure \ref{fig_L}(b)) are exactly the graphs $G$ with $\ld(G) = \lfloor \log n \rfloor = 2$ (resp. $\ltd(G) = \lfloor \log n \rfloor = 2$).

\begin{corollary} \label{cor:LD}
For any integer $k \ge 2$, %$k \ge 1$, 
the total number of graphs $G$ on $2^k-1+k$ vertices and with LD-number $k$ is equal to the total number of graphs on $k$ vertices times the total number of graphs on $2^k-1$ vertices, that is, $\eta(k) \times \eta(2^k-1)$.
\end{corollary}

\begin{proof}
By the characterization in Theorem~\ref{pro_L}, the only graphs on $2^k-1+k$ vertices and LD-number $k$ are all the graphs $G^L(k)$. Therefore, the result follows by their count in Lemma~\ref{rem:GL(k)-cardinality}.
\end{proof}

Corollary~\ref{cor:LD} is a small correction to the result in~\cite{S_1988} where it was mentioned that the total number of graphs on $2^k-1+k$ vertices with LD-number $k$ is the number of graphs on $k$ vertices.

\subsection{The case of OD- and OTD-codes}
\label{Sub_2_O}

\begin{figure}
\begin{center}
\includegraphics[scale=1.0]{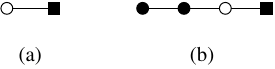}
\caption{Two OD-admissible graphs $K_2$ in (a) and $P_4$ in (b). Black-filled vertices belong to a minimum OD-code, the vertex with empty signature is indicated as square.}
\label{fig_OD-case}
\end{center}
\end{figure}

Here, we construct graphs minimizing their X-density for
$X \in \{OD, OTD\}$. Recall that for open-separating sets $C$, it is required that the open signatures $N(v)\cap C$ are unique sets for all vertices $v$ of the graph. Note that in the case of OD-codes, the signature $N(v)\cap C$ of at most one vertex $v \in C$ can be empty \cite{CW_ISCO2024}; see Figure \ref{fig_OD-case} for two examples of graphs where all minimum OD-codes require a vertex with empty signature. The next two lemmas regarding OD- and OTD-codes appear in~\cite{CW_ISCO2024, CW_ODarxiv2024} and~\cite{SS_2010}, respectively. Nevertheless, as before, we prove them here for the continuity of our arguments and the self-containment of the article.

\begin{lemma}\label{rem:OD-cardinality}
If $G$ is an OD-admissible graph on $n$ vertices and $\od(G) = k$, then $n \le 2^k$. In particular, therefore, we have $\od(G) \ge \lceil \log n \rceil$.
\end{lemma}

\begin{proof}
Let $C$ be a minimum OD-code of $G$. Then, each vertex $v \in V(G) \setminus C$ has a unique non-empty open signature with respect to $C$. Such signatures can be at most $2^{|C|}-1$ in number. In addition, for $C$ to open-separate $v$ from any vertex in $C$, the signature $N(v) \cap C$ must be distinct from any set in $\calN(C) \setminus \{\emptyset\}$. Since $C$ is an open-separating set, any two sets in $\calN(C)$ are distinct from one another. Therefore, we have $|\calN(C) \setminus \{ \emptyset \}| \ge |C|-1$. This implies that $|V(G) \setminus C| \le 2^{|C|}-1 - |\calN(C) \setminus \{ \emptyset \}| \le 2^{|C|} - |C|$. Hence, $n = |V(G)| = |C| + |V(G) \setminus C| = 2^{|C|} \le 2^k$.

Hence, $k \ge \log n$, that is, $k \ge \lceil \log n \rceil$.
\end{proof}

\begin{lemma}\label{rem:OTD-cardinality}
If $G$ is an OTD-admissible graph on $n$ vertices and $\otd(G) = k$, then $n \le 2^k-1$. In particular, therefore, we have $\otd(G) \ge \lceil \log (n+1) \rceil$.
\end{lemma}

\begin{proof}
Let $C$ be a minimum OTD-code of $G$. Then, each vertex $v \in V(G) \setminus C$ has a unique non-empty open signature with respect to $C$. Such signatures can be at most $2^{|C|}-1$ in number. Moreover, for $C$ to open-separate $v$ from any vertex in $C$, the signature $N(v) \cap C$ must be distinct from any set in $\calN(C) \setminus \{\emptyset\}$. In addition, since $C$ is a total-dominating set of $G$, we have $\emptyset \notin \calN_C$ and hence, $\calN(C) \setminus \{ \emptyset \} = \calN(C)$. Furthermore, since $C$ is an open-separating set, any two sets in $\calN(C)$ are distinct from one another. Therefore, we have $|\calN(C) \setminus \{ \emptyset \}| = |\calN(C)| = |C|$. This implies that $|V(G) \setminus C| \le 2^{|C|}-1 - |\calN_C \setminus \{ \emptyset \}| = 2^{|C|} - 1 - |C|$. Therefore, $n = |V(G)| = |C| + |V(G) \setminus C| \le 2^{|C|}-1 = 2^k-1$.

Hence, $k \ge \log (n+1)$, that is, $k \ge \lceil \log (n+1) \rceil$.
\end{proof}

Next, we construct graphs realizing the minimum X-densities for $X \in \{OD, OTD\}$. Define the graph $G^O(k) = (C \cup V \setminus C, E)$ for $k \geq 2$ as follows.
\begin{itemize}
\itemsep0pt
  \item $C$ is a vertex subset of cardinality $k$ and induces a subgraph $G^O(k)[C]$ without open twins,
  \item $V \setminus C$ has exactly one vertex $v_{C'}$ with signature $N(v_{C'}) \cap C = C'$ for each non-empty subset $C' \subset C$ such that there is no vertex $u \in C$ with open signature $N(u) \cap C = C'$.
\end{itemize}

\begin{figure}[h]
\begin{center}
\includegraphics[scale=1.0]{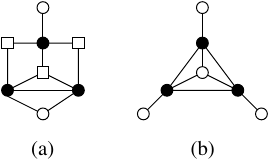}
\caption{The graphs $G^O(3)$. Black-filled vertices belong to the minimum code $C$ (no optional edges between the vertices in $V \setminus C$ are indicated for better readability). The graph in (a) is extremal for OD-codes, the graph in (b) is extremal for OTD-codes.}
\label{fig_O}
\end{center}
\end{figure}

Notice that, on account of being open-twin-free, the subgraph $G^O(k)[C]$ is OD-admissible (refer to Table~\ref{tab_admissible}). 
Also, notice that the construction of $G^O(k)$ imposes conditions only on those edges that are incident with a vertex in $C$, but imposes no conditions on the edge set of $G^O(k)[V \setminus C]$. For examples, the graph $G^O(2)$ is a clique of size 3 and Figure \ref{fig_O} illustrates the graphs $G^O(k)$ for $k = 3$.

\begin{lemma}\label{rem:GO(k)-cardinality}
For the graph $G^O(k) = (C \cup V \setminus C, E)$, we have $|V| = 2^k$ if $G^O(k)[C]$ has an isolated vertex; and $|V| = 2^k-1$ if $G^O(k)[C]$ has no isolated vertex. In particular, the number of graphs $G^O(k)$ is $\overline{\eta}^O(k) \times \eta(2^k-1-k) + \overline{\eta}^O(k-1) \times \eta(2^k-k)$.
\end{lemma}

\begin{proof}
Let us first assume that $G^O(k)[C]$ has an isolated vertex. By construction of $G^O(k)$, since each subset of $C \setminus (\calN(C) \cup \{ \emptyset \})$ is a signature of a vertex in $V \setminus C$, we have $|V \setminus C| = 2^{|C|} - |\calN(C) \cup \{ \emptyset \}|$. Since $G^O(k)[C]$ has an isolated vertex, it implies that $\emptyset \in \calN(C)$ and hence, $\calN(C) \cup \{ \emptyset \} = \calN(C)$. Moreover, $|\calN(C)| = |C|$ since any two sets in $\calN(C)$ are distinct on account of $C$ being an open-separating set of $G^O(k)$ by Lemma~\ref{rem:OD-cardinality}. Hence, $|V \setminus C| = 2^{|C|}-|C|$ and thus, $|V| = |C| + |V \setminus C| = 2^{|C|} = 2^k$.

Let us now assume that $G^O(k)[C]$ has no isolated vertex. In this case, we have $\emptyset \notin \calN(C)$ and hence, $|\calN(C) \cup \{ \emptyset \}| = |\calN(C)|+1$. Then, by the exact same argument as before, we have $|\calN(C)| = |C|$ and $|V \setminus C| = 2^{|C|} - |\calN(C) \cup \{ \emptyset \}| = 2^{|C|}-1 - |\calN(C)| = 2^{|C|}-1-|C|$. Hence, $|V| = |C| + |V \setminus C| = 2^{|C|}-1 = 2^k-1$.

By the first assertion, since $|V \setminus C|=2^k-1-k$ in the case that $G^O(k)[C]$ is isolate-free, the number of graphs $G^O(k)$ with an isolate-free subgraph $G^O(k)[C]$ is $\overline{\eta}^O(k) \times \eta(2^k-1-k)$. Similarly, since $|V \setminus C|=2^k-k$ in the case that $G^O(k)[C]$ has an isolated vertex, the number of graphs $G^O(k)$ with its subgraph $G^O(k)[C]$ having an isolated vertex is $\overline{\eta}^O(k-1) \times \eta(2^k-k)$. This implies the result on the number of graphs $G^O(k)$.
\end{proof}

We next show that, for $X \in \{OD, OTD\}$, the graphs $G^O(k)$ are the ones with the minimum possible X-densities. As it can be easily verified that an X-admissible graph $G$ satisfies $\x(G) = 2$ if and only if $G \cong G^O(2) \cong K_3$, we let $k \ge 3$ in the next theorems of this section.

\begin{theorem}\label{thm_O}
Let $k \ge 3$. Then in $G^O(k) = (C \cup V \setminus C, E)$, $C$ is a minimum%\\[-6mm]
\begin{itemize}
\itemsep0pt
  \item OD-code,
  \item OTD-code if $G^O(k)[C]$ has no isolated vertex, that is, if $G^O(k)[C]$ is OTD-admissible.%\\[-6mm]
\end{itemize}
Moreover, for $X \in \{OD, OTD\}$, the minimum X-density of a graph with X-number $k$ is attained by the graphs $G^O(k)$ with $G^O(k)[C]$ having an isolated vertex for $X=OD$ and by those $G^O(k)$ with isolate-free $G^O(k)[C]$ for $X=OTD$.
\end{theorem}

\begin{proof}
The open-twin-free $G^O(k)[C]$ condition in the construction of $G^O(k)$ implies that each vertex in $C$ has a unique (possibly empty) open signature with respect to $C$. In addition, the second condition makes sure that the open signature of a vertex in $V \setminus C$ with respect to $C$ is different from that of any other vertex either in $C$ or in $V \setminus C$. Therefore, $C$ is an open-separating set of $G^O(k)$. Moreover, since no vertex in $V \setminus C$ has an empty open signature in $C$, the set $C$ total-dominates all vertices of $V \setminus C$. This implies that $C$ is a dominating set (and hence, an OD-code) of $G^O(k)$ and is a total-dominating set (and hence, and OTD-code) of $G^O(k)$ if $G^O(k)[C]$ has no isolated vertices. We now prove that the set $C$ is also of minimum cardinality in both cases. 
If to the contrary, the set $C$ is not of minimum cardinality, then by Lemmas~\ref{rem:OD-cardinality} and~\ref{rem:OTD-cardinality}, we would have $|V| \le 2^{k-1}$ contradicting $|V| \ge 2^k-1$ (by Lemma~\ref{rem:GO(k)-cardinality} and the fact that $2^{k-1} < 2^k-1$, since $k \ge 3$). 

Therefore, we have $\x(G^O(k)) = k$ both for $X=OD$ and for $X=OTD$ with isolate-free $G^O(k)[C]$. For the case that $X=OD$ and $G^O(k)[C]$ has an isolated vertex, we have $|V| = 2^k$ (by Lemma~\ref{rem:GO(k)-cardinality}) and hence, the result that such graphs $G^O(k)$ have the minimum OD-densities among all graphs with OD-number $k$ follows by Lemma~\ref{rem:OD-cardinality}. 
On the other hand, for $X=OTD$ and isolate-free $G^O(k)[C]$, we have $|V| = 2^k-1$ (again by Lemma~\ref{rem:GO(k)-cardinality}) and hence, the result that such graphs $G^O(k)$ have the minimum OTD-densities among all graphs with OTD-number $k$ follows by Lemma~\ref{rem:OTD-cardinality}.
\end{proof}

Next, for $X \in \{OD, OTD\}$, we characterize all X-admissible graphs which attain the logarithmic lower bounds in Lemmas~\ref{rem:OD-cardinality} and~\ref{rem:OTD-cardinality}.

\begin{theorem} \label{pro_OD}
Let $G$ be an OD-admissible graph of order $n$ and let $k = \lceil \log n \rceil \ge 3$. Then we have $\od(G) = k$ if and only if $G$ can be obtained from a $G^O(k) = (C \cup V \setminus C, E)$ by removing up to $2^{k-1}-1$ vertices of $V \setminus C$ if $G^O(k)[C]$ has an isolated vertex and up to $2^{k-1}-2$ vertices of $V \setminus C$ otherwise.
\end{theorem}

\begin{proof}
To prove the sufficiency part, let $G \cong G^O(k) - S$ for some graph $G^O(k)$, where $S \subset V \setminus C$ with  $|S| \le 2^{k-1}-1$ if $G^O(k)[C]$ has an isolated vertex and $|S| \le 2^{k-1}-2$ otherwise. Then, by Theorem~\ref{thm_O}, the set $C$ is a minimum OD-code of $G^O(k)$. Thus, $C$ is also an OD-code of $G^O(k)-S$. Let $D \subset V(G)$ be the isomorphic copy of $C$ in $G$. Then, by the isomorphism, $D$ is an OD-code of $G$. Moreover, we have $|D| = |C| = k = \lceil \log n \rceil$. Hence, by Lemma~\ref{rem:OD-cardinality}, we have $\od(G) = k$. 

To prove the necessary part, let $\od(G) = k$ and let $C$ be a minimum OD-code of $G$. It is clear that every vertex in $V(G) \setminus C$ has a set in $2^C \setminus (\calN(C) \cup \{\emptyset\})$ as its unique non-empty open signature with respect to $C$. However, since the vertices in $V(G) \setminus C$ may not correspond to all the possible sets in $2^C \setminus (\calN(C) \cup \{\emptyset\})$ as their signatures, we therefore have $G \cong G^O(k) - S$, where $S \subset V \setminus C$ and the edge set of $G^O(k)[V \setminus C]$ is so chosen to allow for the isomorphism. Now, if $G^O(k)[C]$ contains (resp. does not contain) an isolated vertex, then by Lemma~\ref{rem:GO(k)-cardinality}, we have $|V| = 2^k$ (resp. $|V| = 2^k-1$). We now show that $|S| \le 2^{k-1}-1$ (resp. $|S| \le 2^{k-1}-2$). Let us now assume to the contrary that $|S| \ge 2^{k-1}$ (resp. $|S| \ge 2^{k-1}-1$). This implies that $|V|-|S| \le 2^{k-1}$ . Since $n = |V| - |S|$ by the isomorphism $G \cong G^O(k) - S$, we have $\log n = \log (|V|-|S|) \le k-1$, that is, $\lceil \log n \rceil < k$, a contradiction to $k = \lceil \log n \rceil$. This proves the necessary part.
\end{proof}

For instance, the graph $G^O(3)$ from Figure \ref{fig_O}(a) (resp. Figure \ref{fig_O}(b)) and all graphs obtained by removing up to 3 (resp. 2) un-filled vertices are exactly the graphs $G$ with $\gamma^{OD}(G) = \lceil \log n \rceil = 3$.

{\begin{remark}\label{rem_OD-case} \rm
As pointed out in \cite{CW_ISCO2024}, a disconnected OD-admissible graph can have an OD-number strictly larger than the sum of the OD-numbers of its components, namely, when the minimum OD-codes of at least two components require a vertex with empty signature. For instance, we have $\od(K_2+K_2)=3$, $\od(K_2+P_4)=5$, $\od(P_4+P_4)=7$. On the contrary, it is well-known that the OTD-number of a graph always equals the sum of the OTD-numbers of its components.

  In the characterization of an extremal graph $G$ in Theorem~\ref{pro_OD}, in order to obtain a disconnected graph by removing vertices from $V \setminus C$ in $G^O(k)$, we note that $G^O(k)[C]$ needs to be disconnected, say it has two components $C_1$ and $C_2$ of cardinality $|C_1|= \ell$ and $|C_2|= k-\ell$.
  To obtain a disconnected graph, all vertices $v_{C'} \in V \setminus C$ must be removed whose signatures $C'$ satisfy $C' \cap C_1 \neq \emptyset$ and $C' \cap C_2 \neq \emptyset$, hence in total
  $(2^{\ell} - 1)(2^{k-\ell} - 1)$ many vertices.
  For $\ell = 1$, we have
  $$
  (2^{1} - 1)(2^{k-1} - 1) = 2 \cdot 2^{k-1} - 2 - 2^{k-1} + 1 = 2^{k-1} - 1
  $$
  such vertices, that is the maximum number of vertices that can be removed by Theorem~\ref{pro_OD} (for all other values of $\ell$, this limit is exceeded).
  Hence, the only case when removing up to $2^{k-1}-1$ vertices from $V \setminus C$ in $G^O(k)$ can result in a disconnected graph $G$ occurs when $G^O(k)[C]$ has one isolated vertex, say $u$, and all $2^{k-1} - 1$ vertices $v_{C' \cup \{u\}} \in V \setminus C$ are removed for all non-empty subsets $C' \subseteq C \setminus \{u\}$.
  Then the resulting graph $G$ has two components: a $K_2$ composed of $u$ and $v_{\{u\}}$, and the graph induced by all vertices in $C \setminus \{u\}$ and all remaining vertices in $(V \setminus C) \setminus \{v_{\{u\}}\}$. In this case, $G$ is isomorphic to $K_2 + G^O(k-1)$ and $\od(G) = 1 + (k-1) = k$ clearly follows.
For instance, the graph obtained from $G^O(3)$ in Figure \ref{fig_O}(a) by removing the 3 unfilled vertices marked as squares is isomorphic to $K_2 + G^O(2)$ and satisfies $\od(G) = 3$.
\end{remark}

\begin{theorem}\label{pro_OTD}
Let $G$ be an OTD-admissible graph of order $n$ and let $k = \lceil \log (n+1) \rceil \dc{\ge 3}$. Then we have $\otd(G) = k$ if and only if $G$ can be obtained from a $G^O(k) = (C \cup V \setminus C, E)$ with an OTD-admissible $G^O(k)[C]$ by removing up to $2^{k-1}-1$ vertices of $V \setminus C$.
\end{theorem}

\begin{proof}
To prove the sufficiency part, let $G \cong G^O(k) - S$ for some graph $G^O(k)$, where the subgraph $G^O(k)[C]$ is OTD-admissible and $S \subset V \setminus C$ with  $|S| \le 2^{k-1}-1$. Then, by Theorem~\ref{thm_O}, the set $C$ is a minimum OTD-code of $G^O(k)$. Thus, $C$ is also an OTD-code of $G^O(k)-S$. Let $D \subset V(G)$ be the isomorphic copy of $C$ in $G$. Then, by the isomorphism, $D$ is an OTD-code of $G$. Moreover, we have $|D| = |C| = k = \lceil \log (n+1) \rceil$. Hence, by Lemma~\ref{rem:OTD-cardinality}, we have $\otd(G) = k$.

To prove the necessary part, let $\otd(G) = k$ and let $C$ be any minimum OTD-code of $G$. Then,  every vertex in $V(G) \setminus C$ has a set in $2^C \setminus (\calN(C) \cup \{\emptyset\})$ as its unique non-empty open signature with respect to $C$. However, since the vertices in $V(G) \setminus C$ may not correspond to all the possible sets in $2^C \setminus (\calN(C) \cup \{\emptyset\})$ as their signatures, we therefore have $G \cong G^O(k) - S$, where $S \subset V \setminus C$ and the edge set of $G^O(k)[V \setminus C]$ is so chosen to allow for the isomorphism. Also, $G^O(k)[C]$ is OTD-admissible since $C$ is an OTD-code of $G$. Therefore, by Lemma~\ref{rem:GO(k)-cardinality}, we have $|V| = 2^k-1$. We now show that $|S| \le 2^{k-1}-1$. Let us assume to the contrary that $|S| \ge 2^{k-1}$. This implies that $|V|-|S| \le 2^{k-1}-1$. Since $n = |V| - |S|$ by the isomorphism $G \cong G^O(k) - S$, we have $\log (n+1) = \log (|V|-|S|+1) \le k-1$, that is, $\lceil \log (n+1) \rceil < k$, a contradiction to $k = \lceil \log (n+1) \rceil$. This proves the necessary part.
\end{proof}

For instance,
the graph $G^O(3)$ from Figure \ref{fig_I}(b) and all graphs obtained by removing up to 2 unfilled vertices are exactly the graphs $G$ with $\gamma^{OTD}(G) = \lceil \log (n+1) \rceil = 3$.

\subsection{The case of ID- and ITD-codes}
\label{Sub_2_I}

In this subsection, we construct graphs minimizing their X-density for $X \in \{ID, ITD\}$. Recall that for closed-separating sets $C$, it is required that the closed signatures $N[v]\cap C$ are unique sets for all vertices $v$ of the graph.

\begin{lemma}\label{rem:I-cardinality}
For $X \in \{ID, ITD\}$, if $G$ is an X-admissible graph on $n$ vertices and with $\x(G) = k$, then $n \le 2^k-1$. In particular, therefore, we have $\x(G) \ge \lceil \log (n+1) \rceil$.
\end{lemma}

\begin{proof}
The result follows by the fact that, if $C$ is a minimum X-code of $G$, then there can be at most $2^{|C|}-1$ non-empty subsets of $C$.
\end{proof}

The result in Lemma~\ref{rem:I-cardinality} for $X=ID$ also appears in~\cite{KCL_1998}. Moreover, Lemma~\ref{rem:I-cardinality} also establishes the corresponding logarithmic lower bound for ITD-numbers of graphs. Next, we construct graphs realizing the minimum X-densities for $X \in \{ID, ITD\}$.
Define the graph $G^I(k) = (C \cup V \setminus C, E)$ for $k \geq 2$ as follows.
\begin{itemize}
\itemsep0pt
  \item $C$ is a vertex subset of cardinality $k$ and induces a subgraph $G^I(k)[C]$ without closed twins,
  \item $V \setminus C$ has exactly one vertex $v_{C'}$ with signature $N(v_{C'}) \cap C = C'$ for each non-empty subset $C' \subseteq C$ such that there is no vertex $u \in C$ with closed signature $N[u]\cap C = C'$.
\end{itemize}

Notice that, on account of being closed-twin-free, the subgraph $G^I(k)[C]$ is ID-admissible (refer to Table~\ref{tab_admissible}). 
Also, as in the previous cases, notice that the construction of $G^I(k)$ imposes conditions only on those edges that are incident with a vertex in $C$, but imposes no condition on the edge set of $G^I(k)[V \setminus C]$. For instance, $G^I(2)$ is a path on three vertices and Figure \ref{fig_I} depicts the possible graphs $G^I(k)$ for $k = 3$. 

\begin{lemma} \label{rem:GI(k)-cardinality}
For the graph $G^I(k) = (C \cup V \setminus C, E)$, we have $|V| = 2^k-1$. In particular, the number of graphs $G^I(k)$ is $\eta^I(k) \times \eta(2^k-1-k)$.
\end{lemma}

\begin{proof}
By the construction of $G^I(k)$, for each non-empty subset $C'$ of $C$, the graph $G^I(k)$ has (exactly) one vertex with $C'$ as it closed signature. Therefore, the first assertion follows by the fact that there can be at most $2^{|C|}-1 = 2^k-1$ non-empty such subsets of $C$.

Moreover, the second assertion follows by the fact that the graphs $G^I(k)[C]$, admitting $I$-separating sets, are 
$\eta^I(k)$ %$\eta(k)$ 
in number; and the graphs $G^I(k)[V \setminus C]$, with no restrictions on their edge sets, are $\eta(2^k-1-k)$ in number.
\end{proof} 

We next show that, for $X \in \{ID, ITD\}$, the graphs $G^I(k)$ are the ones with the minimum possible X-densities. 
It can be easily verified that there is no ITD-admissible graph with ITD-number~$2$ and that the only ID-admissible graph $G$ with ID-number~$2$ is $G \cong G^I(2) \cong P_3$. 
Therefore, for $X \in \{ID, ITD\}$, we let $k \ge 3$ in the next two theorems.

\begin{figure}
\begin{center}
\includegraphics[scale=1.0]{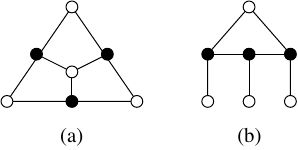}
\caption{The two graphs $G^I(3)$. Black-filled vertices belong to the minimum code $C$ (no optional edges between the vertices in $V \setminus C$ are indicated for better readability). The graph in (a) is extremal for ID-codes, the graph in (b) is extremal for both ID- and ITD-codes.}
\label{fig_I}
\end{center}
\end{figure}

\begin{theorem}\label{thm_I}
Let $k \ge 3$. Then in $G^I(k) = (C \cup V \setminus C, E)$, $C$ is a minimum
\begin{itemize}
\itemsep0pt
  \item ID-code,
  \item ITD-code if $G^I(k)[C]$ has no isolated vertex, that is, if $G^I(k)[C]$ is ITD-admissible.
\end{itemize}
Moreover, for $X \in \{ID, ITD\}$, the minimum X-density of a graph with X-number $k$ is attained by any $G^I(k)$ for $X=ID$ and by those $G^I(k)$ with isolate-free $G^I(k)[C]$ for $X=ITD$.
\end{theorem}

\begin{proof}
The closed-twin-free $G^I(k)[C]$ condition in the construction of $G^I(k)$ implies that each vertex in $C$ has a unique closed signature with respect to $C$. In addition, the second condition makes sure that the closed signature of a vertex in $V \setminus C$ with respect to $C$ is different from that of any other vertex either in $C$ or in $V \setminus C$. Therefore, $C$ is a closed-separating set of $G^I(k)$. Moreover, since no vertex in $V \setminus C$ has an empty open signature in $C$, the set $C$ total-dominates all vertices of $V \setminus C$. This implies that $C$ is a dominating set (and hence, an ID-code) of $G^I(k)$ and is a total-dominating set (and hence, and ITD-code) of $G^I(k)$ if $G^I(k)[C]$ has no isolated vertices. Moreover, in both cases, the code $C$ must be of minimum cardinality, as otherwise, by Lemma~\ref{rem:I-cardinality}, we would have $|V| \le 2^{k-1}-1$ contradicting $|V| = 2^k-1$ (by Lemma~\ref{rem:GI(k)-cardinality} and the fact that $2^{k-1} < 2^k$ for all $k \ge 3$).

Therefore, we have $\x(G^I(k)) = k$ both for $X=ID$ and for $X=ITD$, with isolate-free $G^I(k)[C]$. Thus, for $X \in \{ID, ITD\}$ the result that such graphs $G^I(k)$ have the minimum X-densities among all graphs with X-number $k$ follows  by the comparison of graph orders in Lemmas~\ref{rem:I-cardinality} and~\ref{rem:GI(k)-cardinality}.
\end{proof}

Next, for $X \in \{ID, ITD\}$, we characterize all X-admissible graphs which attain the logarithmic lower bound in Theorem~\ref{rem:I-cardinality}. Even though the same characterization has also appeared for $X=ID$ in~\cite{M_2006}, we still include a proof of the result as part of the general construction technique undertaken in this article.

\begin{theorem}\label{pro_I}
For $X \in \{ID, ITD\}$, let $G$ be an X-admissible graph of order $n$ and let $k = \lceil \log (n+1) \rceil \ge 3$. Then we have $\x(G) = k$ if and only if $G$ can be obtained from a $G^I(k) = (C \cup V\setminus C, E)$ with an X-admissible $G^I(k)[C]$ by removing at most $2^{k-1}-1$ vertices of $V \setminus C$.
\end{theorem}

\begin{proof}
To prove the sufficiency part, let $G \cong G^I(k) - S$ for some $G^I(k)$, where  $G^I(k)[C]$ is X-admissible and $S \subset V \setminus C$ with $|S| \le 2^{k-1}-1$. Then, by Theorem~\ref{thm_I}, the set $C$ is a minimum X-code of $G^I(k)$. Thus, $C$ is also an X-code of $G^I(k)-S$. Let $D \subset V(G)$ be the isomorphic copy of $C$ in $G$. Then, by the isomorphism, $D$ is an X-code of $G$. Moreover, we have $|D| = |C| = k = \lceil \log (n+1) \rceil$. Hence, by Lemma~\ref{rem:I-cardinality}, we have $\x(G) = k$.

To prove the necessary part, let us assume that $\x(G) = k$ and let $C$ be any minimum X-code of $G$. Then, it is clear that every vertex in $V(G) \setminus C$ has a set in $2^C \setminus (\calN[C] \cup \{ \emptyset \})$ as its unique non-empty closed signature with respect to $C$. However, since the vertices in $V(G) \setminus C$ may not correspond to all the possible sets in $2^C \setminus (\calN[C] \cup \{ \emptyset \})$ as their signatures, we therefore have  $G \cong G^I(k) - S$, where $S \subset V \setminus C$ and the edge set of $G^I(k)[V \setminus C]$ is so chosen to allow for the isomorphism. Also, $G^I(k)[C]$ is X-admissible since $C$ is an X-code of $G$. We now show that $|S| \le 2^{k-1}-1$. Let us assume to the contrary that $|S| \ge 2^{k-1}$. By Lemma~\ref{rem:GI(k)-cardinality}, we have $|V|= 2^k-1$. Therefore, $|V|-|S| \le 2^{k-1}-1$. Since $n = |V| - |S|$ by the isomorphism $G \cong G^I(k) - S$, we have $\lceil \log (n+1) \rceil = \lceil \log (|V|-|S|+1) \rceil \le k-1$, a contradiction to $k = \lceil \log (n+1) \rceil$. This proves the necessary part.
\end{proof}

For instance, the graph $G^I(3)$ and all graphs obtained by removing up to three unfilled vertices from $G^I(3)$ in Figure \ref{fig_I} (resp. Figure \ref{fig_I}(b)) are exactly the graphs $G$ with $\id(G) = \lceil \log (n+1) \rceil
= 3$ (resp. $\itd(G) = \lceil \log (n+1) \rceil = 3$).

\subsection{The case of FD- and FTD-codes}
\label{Sub_2_F}

In this section, we construct graphs minimizing their X-density for $X \in \{FD, FTD\}$. Note that in the case of FD-codes, the open signature $N(v)\cap C$ of one vertex $v \in C$ can be empty \cite{CW_2024-FS}, see Figure \ref{fig_FD-case} for an example of a graph such that all minimum FD-codes require a vertex with empty signature. As further pointed out in \cite{CW_2024-FS}, a disconnected FD-admissible graph can have an FD-number strictly larger than the sum of the FD-numbers of its components, namely, when the minimum FD-codes of at least two components require a vertex with empty signature. For instance, we have $\fd(B_3+B_3)=11$. On the contrary, for a disconnected FTD-admissible graph, its FTD-number always equals the sum of the FTD-numbers of its components by \cite{CW_2024-FS}.

Recall that for full-separating sets $C$, it is required that $N(v)\cap C \neq N(u)\cap C$ and $N[v]\cap C \neq N[u]\cap C$ for all distinct vertices $v,u$ of the graph. These requirements imply the following series of results.

\begin{lemma}\label{rem:full-sep}
Let $C$ be a full-separating set of a graph $G$. Then, we have $\calN(C) \cap \calN[C] = \emptyset$. Moreover, we have $|\calN(C)| = |\calN[C]| = |C|$. In particular, therefore, $|\calN_C| = 2|C|$. 
\end{lemma}

\begin{proof}
Let us assume to the contrary that $\calN(C) \cap \calN[C] \ne \emptyset$. This implies that $N(u) \cap C = N[v] \cap C$ for some distinct vertices $u,v \in C$. In other words, $N[v] \cap C \subseteq N[u] \cap C$. Moreover, this also implies that $uv$ is an edge in $G$ and hence, $u \in N[v]$. Therefore, we have $N[u] \cap C \subseteq N[v] \cap C$ as well, that is, $N[u] \cap C = N[v] \cap C$ contradicting the fact that $C$ full-separates the pair $u,v$. This proves the first result.

Since, $C$ is a full-separating set of $G$, each pair of sets in either $\calN(C)$ or $\calN[C]$ must be distinct. Therefore, $|\calN(C)| = |\calN[C]| = |C|$. Thus, using $\calN(C) \cap \calN[C] = \emptyset$, we obtain $|\calN_C| = 2|C|$. 
\end{proof}

\begin{lemma}\label{rem:FD-cardinality}
If $G$ is an FD-admissible graph on $n$ vertices and $\fd(G) = k$, then $n \le 2^k - k$. In particular, therefore, we have $\fd(G) \ge 1 + \lfloor \log n \rfloor$.
\end{lemma}

\begin{proof}
Let $C$ be a minimum FD-code of $G$. Then, each vertex $v \in V(G) \setminus C$ has a unique non-empty open signature with respect to $C$. Such signatures can be at most $2^{|C|}-1$ in number. In addition, for $C$ to full-separate $v$ from any vertex in $C$, the signature $N(v) \cap C$ must be distinct from any set in $\calN_C \setminus \{\emptyset\}$. By Lemma~\ref{rem:full-sep}, we have $|\calN_C \setminus \{ \emptyset \}| \ge 2|C|-1$. This implies that $|V(G) \setminus C| \le 2^{|C|}-1 - |\calN_C \setminus \{ \emptyset \}| \le 2^{|C|} - 2|C| = 2^k-2k$. Hence, $n = |V(G)| = |C| + |V(G) \setminus C| \le 2^k-k$.

Since $2^k - k < 2^k$ for all $k \ge 1$, we have $k > \log n$, that is, $k \ge 1+\lfloor \log n \rfloor$.
\end{proof}

\begin{lemma}\label{rem:FTD-cardinality}
If $G$ is an FTD-admissible graph on $n$ vertices and $\ftd(G) = k$, then $n \le 2^k-1-k$. In particular, therefore, we have $\ftd(G) \ge 1 + \lfloor \log (n+1) \rfloor$.
\end{lemma}

\begin{proof}
Let $C$ be a minimum FTD-code of $G$. Then, each vertex $v \in V(G) \setminus C$ has a unique non-empty open signature with respect to $C$. In addition, for $C$ to full-separate $v$ from any vertex in $C$, the signature $N(v) \cap C$ must be distinct from any set in $\calN_C \setminus \{\emptyset\}$. Since $C$ is a total-dominating set of $G$, we have $\emptyset \notin \calN_C$ and hence, $\calN_C \setminus \{ \emptyset \} = \calN_C$. By Lemma~\ref{rem:full-sep} therefore, we have $|\calN_C \setminus \{ \emptyset \}| = |\calN_C| = 2|C|$. This implies that $|V(G) \setminus C| \le 2^{|C|}-1 - |\calN_C \setminus \{ \emptyset \}| = 2^{|C|} - 1 - 2|C| = 2^k-1-2k$. Hence, $n = |V(G)| = |C| + |V(G) \setminus C| \le 2^k-1-k$.

Since $2^k-1-k < 2^k-1$ for all $k \ge 1$, we have $k > \log (n+1)$, that is, $k \ge 1+\lfloor \log (n+1) \rfloor$.
\end{proof}

\begin{figure}
\begin{center}
\includegraphics[scale=1.0]{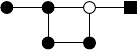}
\caption{The FD-admissible graph $B_3$. Black-filled vertices belong to a minimum FD-code, the vertex with empty open signature is indicated as square.}
\label{fig_FD-case}
\end{center}
\end{figure}

For $X \in \{FD, FTD\}$, Lemmas~\ref{rem:FD-cardinality} and~\ref{rem:FTD-cardinality} establish the logarithmic lower bounds for the X-numbers of X-admissible graphs. Next, we construct graphs realizing the minimum X-densities. Recall from Remark \ref{rem:G[C] X-admissible} that, for any X-code $C$ of an X-admissible graph $G$,
the subgraph $G[C]$ is itself X-admissible. 
For $X \in \{FD, FTD\}$, this means that $G[C]$ is twin-free. As there are no twin-free graphs of order 2 and 3, $G[C]$ needs to have at least 4 vertices. Define the graph $G^F(k) = (C \cup V \setminus C, E)$ for $k \geq 4$ as follows.
\begin{itemize}
\itemsep0pt
  \item $C$ is a vertex subset of cardinality $k$ and induces a subgraph $G^F(k)[C]$ having neither open nor closed twins,
  \item $V \setminus C$ has exactly one vertex $v_{C'}$ with signature $N(v_{C'}) \cap C = C'$ for each non-empty subset $C' \subset C$ such that there is no vertex $u \in C$ with either $N(u) \cap C = C'$ or $N[u]\cap C = C'$.
\end{itemize}

On account of being twin-free, the subgraph $G^F(k)[C]$ is indeed FD-admissible (refer to Table~\ref{tab_admissible}). Moreover, as in the previous cases, notice that the construction of $G^F(k)$ imposes conditions only on those edges that are incident with a vertex in $C$, but  %; and 
imposes no condition on the edge set of $G^F(k)[V \setminus C]$. For illustration, Figure \ref{fig_F} shows the possible graphs $G^F(k)$ for $k = 4$ (with $G^F(k)[C]$ isomorphic to the $P_4$ as only twin-free graph on 4 vertices).

\begin{lemma}\label{rem:GF(k)-cardinality}
For the graph $G^F(k) = (C \cup V \setminus C, E)$, we have $|V| = 2^k-k$ if $G^F(k)[C]$ has an isolated vertex; and $|V| = 2^k-1-k$ if $G^F(k)[C]$ has no isolated vertex. In particular, the number of graphs $G^F(k)$ is $\overline{\eta}^F(k) \times \eta(2^k-1-2k) + \overline{\eta}^F(k-1) \times \eta(2^k-2k)$.
\end{lemma}

\begin{proof}
Let us first assume that $G^F(k)[C]$ has an isolated vertex. Since each set in $2^C \setminus (\calN_C \cup \{ \emptyset \})$ is a signature of a vertex in $V \setminus C$, we have $|V \setminus C| = 2^{|C|} - |\calN_C \cup \{ \emptyset \}|$. As $G^F(k)[C]$ has an isolated vertex, it implies that $\emptyset \in \calN(C)$ and hence, $\calN_C \cup \{ \emptyset \} = \calN_C$. Therefore, by Lemma~\ref{rem:full-sep}, we have $|\calN_C \cup \{ \emptyset \}| = |\calN_C| = 2|C|$. Hence, $|V \setminus C| = 2^{|C|}-2|C|$ and thus, $|V| = |C| + |V \setminus C| = 2^{|C|}-|C| = 2^k-k$.

Let us now assume that $G^F(k)[C]$ has no isolated vertex. In this case, we have $\emptyset \notin \calN_C$ and hence, $|\calN_C \cup \{ \emptyset \}| = |\calN_C|+1$. Then, by the exact same argument as before, we have $|V \setminus C| = 2^{|C|} - |\calN_C \cup \{ \emptyset \}| = 2^{|C|}-1-2|C|$ (using Lemma~\ref{rem:full-sep} again). Hence, $|V| = |C| + |V \setminus C| = 2^{|C|}-1-|C| = 2^k-1-k$. 

By the first assertion, since $|V \setminus C|=2^k-1-2k$ in the case that $G^F(k)[C]$ is isolate-free, the number of graphs $G^F(k)$ with an isolate-free subgraph $G^F(k)[C]$ is $\overline{\eta}^F(k) \times \eta(2^k-1-2k)$. Similarly, since $|V \setminus C|=2^k-2k$ in the case that $G^F(k)[C]$ has an isolated vertex, the number of graphs $G^F(k)$ with its subgraph $G^F(k)[C]$ having an isolated vertex is $\overline{\eta}^F(k-1) \times \eta(2^k-2k)$. This implies the result on the number of graphs $G^F(k)$.
\end{proof}

\begin{figure}[h]
\begin{center}
\includegraphics[scale=1.0]{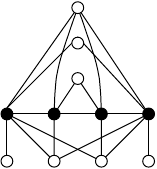}
\caption{The graph $G^F(4)$. Black-filled vertices belong to the minimum code $C$ (no optional edges between the vertices in $V \setminus C$ are indicated for better readability). The graph is extremal for FTD-codes.}
\label{fig_F}
\end{center}
\end{figure}

We next show that, for $X \in \{FD, FTD\}$, the graphs $G^F(k)$ are the ones with the minimum possible X-densities. Since the only twin-free graph (up to isomorphism) on~$4$ vertices is $P_4$, following the construction of $G^F(k)$, it can be verified that for $k=4$, the graphs $G^F(k)$ are determined to be exactly of the form depicted in Figure~\ref{fig_F}, where the subgraph $G^F(k)[V \setminus C]$ can be endowed with any possible edge set. It can also be verified that the only X-admissible graphs $G$ satisfying $\x(G) = 4$ are, up to isomorphism, those that can be obtained by removing up to any~$3$ vertices from $V \setminus C$ of $G^F(4)$ in Figure~\ref{fig_F}. Therefore, we assume $k \ge 5$ in the next theorems of this section.

\begin{theorem}\label{thm_F}
Let $k \ge 5$. Then in $G^F(k) = (C \cup V \setminus C, E)$, $C$ is a minimum
\begin{itemize}
\itemsep0pt
  \item FD-code,
  \item FTD-code if $G^F(k)[C]$ has no isolated vertex, that is, if $G^F(k)[C]$ is FTD-admissible.
\end{itemize}
Moreover, for $X \in \{FD, FTD\}$, the minimum X-density of a graph with X-number $k$ is attained by the graphs $G^F(k)$ with $G^F(k)[C]$ having an isolated vertex for $X=FD$; and by those $G^F(k)$ with isolate-free $G^F(k)[C]$ for $X=FTD$.
\end{theorem}

\begin{proof}
The twin-free $G^F(k)[C]$ condition in the construction of $G^F(k)$ implies that two distinct vertices in $C$ having distinct (possibly empty) open signatures and distinct closed signatures and with respect to $C$. In addition, the second condition makes sure that the open (and closed) signature of a vertex in $V \setminus C$ with respect to $C$ is different from both the open and closed signatures of any other vertex either in $C$ or in $V \setminus C$. Therefore, $C$ is a full-separating set of $G^F(k)$. Moreover, since no vertex in $V \setminus C$ has an empty open signature in $C$, the set $C$ total-dominates all vertices of $V \setminus C$. This implies that $C$ is a dominating set (and hence, and FD-code) of $G^F(k)$ and is a total-dominating set (and hence, and FTD-code) of $G^F(k)$ if $G^F(k)[C]$ has no isolated vertices. We now show that $C$ is also the minimum code in both the cases. If to the contrary, the set $C$ is not of minimum cardinality, then by Lemmas~\ref{rem:FD-cardinality} and~\ref{rem:FTD-cardinality}, we would have $|V| \le 2^{k-1}-(k-1)$ contradicting $|V| \ge 2^k-1-k$ (by Lemma~\ref{rem:GF(k)-cardinality} and the fact that $2^{k-1}-(k-1) < 2^k-1-k$, since $k \ge 5$).

Thus, we have $\x(G^F(k)) = k$ both for $X=FD$ and for $X=FTD$ with isolate-free $G^F(k)[C]$. Then, for the case that $X=FD$ and $G^F(k)[C]$ having an isolated vertex, we have $|V| = 2^k-k$ (by Lemma~\ref{rem:GF(k)-cardinality}) and hence, the result that such graphs $G^F(k)$ have the minimum FD-densities among all graphs with FD-number $k \ne 4$ follows by Lemma~\ref{rem:FD-cardinality}. On the other hand, for the case that $X=FTD$ with isolate-free $G^F(k)[C]$, we have $|V| = 2^k-1-k$ (again by Lemma~\ref{rem:GF(k)-cardinality}) and hence, the result that such graphs $G^F(k)$ have the minimum FTD-densities among all graphs with FTD-number $k$ follows by Lemma~\ref{rem:FTD-cardinality}.
\end{proof}

Next, for $X \in \{FD, FTD\}$, we characterize all X-admissible graphs which attain the logarithmic lower bound in Lemmas~\ref{rem:FD-cardinality} and~\ref{rem:FTD-cardinality}.

\begin{theorem}\label{pro_FD}
Let $G$ be an FD-admissible graph of order $n$ and let $k = 1+\lfloor \log n \rfloor \ge 5$. Then we have $\fd(G) = k$ if and only if $G$ is obtained from $G^F(k) = (C \cup V \setminus C, E)$ by removing up to $2^{k-1}-k$ vertices of $V \setminus C$ if $G^F(k)[C]$ has an isolated vertex; and up to $2^{k-1}-1-k$ vertices of $V \setminus C$ otherwise.
\end{theorem}

\begin{proof}
 To prove the sufficiency part, let $G \cong G^F(k) - S$ for some $G^F(k)$, where $S \subset V \setminus C$ with $|S| \le 2^{k-1}-k$ if $G^F(k)[C]$ has an isolated vertex and $|S| \le 2^{k-1}-1-k$ otherwise. Then, by Theorem~\ref{thm_F}, the set $C$ is a minimum FD-code of $G^F(k)$. Thus, $C$ is also an FD-code of $G^F(k)-S$. Let $D \subset V(G)$ be the isomorphic copy of $C$ in $G$. Then, by the isomorphism, $D$ is an FD-code of $G$. Moreover, we have $|D| = |C| = k = 1 + \lfloor \log n \rfloor$. Hence, by Lemma~\ref{rem:FD-cardinality}, we have $\fd(G) = k$.

To prove the necessary part, let $\fd(G) = k$ and let $C$ be a minimum FD-code of $G$. It is clear that every vertex in $V(G) \setminus C$ has a set in $2^C \setminus (\calN_C \cup \{\emptyset\})$ as its unique non-empty open signature with respect to $C$. However, since the vertices in $V(G) \setminus C$ may not correspond to all the possible sets in $2^C \setminus (\calN_C \cup \{\emptyset\})$ as their signatures, we therefore have $G \cong G^F(k) - S$, where $S \subset V \setminus C$ and the edge set of $G^F(k)[V \setminus C]$ is so chosen to allow for the isomorphism. Now, if $G^F(k)[C]$ contains (resp. does not contain) an isolated vertex, then by Lemma~\ref{rem:GF(k)-cardinality}, we have $|V| = 2^k-k$ (resp. $|V| = 2^k-1-k$). We now show that $|S| \le 2^{k-1}-k$ (resp. $|S| \le 2^{k-1}-1-k$). Let us assume to the contrary that $|S| \ge 2^{k-1}-k+1$ (resp. $|S| \le 2^{k-1}-k$). This implies that $|V|-|S| \le 2^{k-1}-1 < 2^{k-1}$. Since $n = |V| - |S|$ by the isomorphism $G \cong G^F(k) - S$, we have $\log n = \log (|V|-|S|) < k-1$, that is, $1+\lfloor \log n \rfloor < k$, a contradiction to $k = 1+\lfloor \log n \rfloor$. This proves the necessary part.
\end{proof}

{\begin{remark}\rm
  We note that none of the graphs $G$ resulting from $G^F(k)$ by removing up to $2^{k-1}-k$ vertices from $V \setminus C$ can be disconnected 
  for the following reason. At least $2^{k-1}-1$ vertices have to be removed from $V \setminus C$ to obtain a disconnected graph (following the same argumentation from Remark \ref{rem_OD-case}), which is not possible since $2^{k-1}-k < 2^{k-1}-1$ for all $k > 1$.
\end{remark}

\begin{theorem}\label{pro_FTD}
Let $G$ be an FTD-admissible graph of order $n$ and let $k = 1+\lfloor \log (n+1) \rfloor \ge 5$. Then we have $\ftd(G) = k$ if and only if $G$ can be obtained from $G^F(k) = (C \cup V \setminus C,E)$  with an FTD-admissible $G^F(k)[C]$ by removing up to $2^{k-1}-k$ vertices of $V \setminus C$.
\end{theorem}

\begin{proof}
To prove the sufficiency part, let $G \cong G^F(k) - S$ for some $G^F(k)$, where the subgraph $G^F(k)[C]$ of $G^F(k) = (C \cup V \setminus C, E)$ is FTD-admissible and $S \subset V \setminus C$ with  $|S| \le 2^{k-1}-k$. Then, by Theorem~\ref{thm_F}, the set $C$ is a minimum FTD-code of $G^O(k)$. Thus, $C$ is also an FTD-code of $G^F(k)-S$. Let $D \subset V(G)$ be the isomorphic copy of $C$ in $G$. Then, by the isomorphism, $D$ is an FTD-code of $G$. Moreover, we have $|D| = |C| = k = 1 + \lfloor \log (n+1) \rfloor$. Hence, by Lemma~\ref{rem:FTD-cardinality}, we have $\ftd(G) = k$.

To prove the necessary part, let $\ftd(G) = k$ and let $C$ be any minimum FTD-code of $G$. Then, every vertex in $V(G) \setminus C$ has a set in $2^C \setminus (\calN_C \cup \{\emptyset\})$ as its unique non-empty open signature with respect to $C$. However, since the vertices in $V(G) \setminus C$ may not correspond to all the possible sets in $2^C \setminus (\calN_C \cup \{\emptyset\})$ as their signatures, we therefore have $G \cong G^F(k) - S$, where $S \subset V \setminus C$ and the edge set of $G^F(k)[V \setminus C]$ is so chosen to allow for the isomorphism. Also, $G^F(k)[C]$ has no isolated vertex since $C$ is an FTD-code of $G$. Therefore, by Lemma~\ref{rem:GF(k)-cardinality}, we have $|V| = 2^k-1-k$. We now show that $|S| \le 2^{k-1}-k$. Let us assume to the contrary that $|S| \ge 2^{k-1}-k+1$. This implies that $|V|-|S| \le 2^{k-1}-2 < 2^{k-1}-1$. Since $n = |V| - |S|$ by the isomorphism $G \cong G^F(k) - S$, we have $\log (n+1) = \log (|V|-|S|+1) < k-1$, that is, $1 + \lfloor \log (n+1) \rfloor < k$, a contradiction to $k =1 + \lfloor \log (n+1) \rfloor$. This proves the necessary part.
\end{proof}

For instance,
the graph $G^F(4)$ from Figure \ref{fig_F} and all graphs obtained by removing up to 3 unfilled vertices are exactly the graphs $G$ with $\gamma^{FD}(G) = 1 + \lfloor \log n \rfloor = 4$ and $\gamma^{FTD}(G) = 1 + \lfloor \log (n+1) \rfloor = 4$.

\section{Concluding remarks}

Generalizing the results in \cite{M_2006} and~\cite{S_1988} on ID-codes and LD-codes, respectively, we provided in this paper a common construction technique for extremal graphs $G^S(k) = (C \cup V \setminus C, E)$ for the four studied separation properties $S \in \{L,O,I,F\}$ such that $C$ is a minimum X-code of cardinality $k$
and the order $|V|$ of $G^S(k)$ is maximum possible. Thus, summarizing the results obtained in Section \ref{Sec_2}, we have the following corollaries.

\begin{corollary} \label{cor:summary}
Let $X \in \codes$ involving S-separation for some $S \in \{L,O,I,F\}$. If $G$ is an X-admissible graph on $n$ vertices with $\x(G) = k$, then we have
\[
n \le \left\{
\begin{array}{ll}
  2^k+k-1 & \mbox{if } X \in \{LD,LTD\}  \\
  2^k     & \mbox{if } X = OD \\
  2^k-1   & \mbox{if } X \in \{OTD,ID,ITD\} \\
  2^k-k   & \mbox{if } X = FD \\
  2^k-k-1 & \mbox{if } X = FTD \\
\end{array}
\right.
\]
Moreover, the equalities in the above cases are achieved if and only if $G \cong G^S(k) = (C \cup V \setminus C, E)$ with an X-admissible $G^S(k)[C]$ and the latter having an isolated vertex if $S \in \{O,F\}$.
\end{corollary}

\begin{corollary}\label{cor:log-LB}
For $X \in \codes$, let $G$ be an X-admissible graph on $n$ vertices. Then we have
\[
\x(G) \geq \left\{
\begin{array}{ll}
  \lfloor \log n \rfloor        & \mbox{if } X \in \{LD,LTD\}  \\
  \lceil \log n \rceil          & \mbox{if } X = OD \\
  \lceil \log (n+1) \rceil      & \mbox{if } X \in \{OTD,ID,ITD\} \\
  1 + \lfloor \log n \rfloor    & \mbox{if } X = FD \\
  1 + \lfloor \log (n+1) \rfloor& \mbox{if } X = FTD \\
\end{array}
\right.
\]
\end{corollary}}}

For $X \in \{LD,LTD,OD,OTD,ID\}$, the results in Corollary~\ref{cor:summary} already exist in \cite{CW_ISCO2024,HR_2012,KCL_1998,S_1988,SS_2010}. However, Corollary~\ref{cor:summary} newly establishes logarithmic lower bounds for X-numbers with $X \in \{ITD,FD,FTD\}$. We note that these bounds reflect the relations between the X-numbers for all $X \in \codes$ established in \cite{CW_2024-FS}, see Figure \ref{fig_relations}, and the special relations between OD- and OTD-numbers (resp. the FD- and FTD-numbers) found in \cite{CW_ISCO2024} (resp. \cite{CW_2024-FS}). 
In addition, we characterized for all $X \in \codes$ the graphs $G$ such that $\x(G)$ attains the corresponding logarithmic lower bound in Corollary~\ref{cor:summary}.

Following the characterizations of the graphs with X-numbers attaining the logarithmic lower bounds in Corollary~\ref{cor:summary}, Table~\ref{tab_log LB_tight} provides some immediate graph families whose members serve as tight examples for these lower bounds. These graphs from the specific graph families in Table~\ref{tab_log LB_tight} are constructed by choosing the appropriate edge set of $G^S(k)[C]$ and $G^S(k)[V \setminus C]$ (as clique or as independent set) in the characterizing Theorems~\ref{pro_L},~\ref{pro_OD},~\ref{pro_OTD} and~\ref{pro_I}.

A \emph{biparatite graph} is one whose vertex set can be partitioned into two independent sets. A \emph{cobiparatite graph} is one whose vertex set can be partitioned into two cliques. Finally, a \emph{split graph} is one whose vertex set can be partitioned into one independent set and one clique. For $k \ge 2$, the graphs $G^L(k)$ (respectively, $G^I(k)$) with empty edge sets on the subgraphs $G^L(k)[C]$ and $G^L(k)[V \setminus C]$ (respectively, $G^I(k)[C]$ and $G^I(k)[V \setminus C]$) provide tight examples of bipartite graphs whose LD-numbers (respectively, ID-numbers) attain the logarithmic lower bound in Corollary~\ref{cor:log-LB}. Similarly, for $k \ge 2$, the graphs $G^L(k)$ (respectively, $G^O(k)$) with a complete edge set on the subgraphs $G^L(k)[C]$ and $G^L(k)[V \setminus C]$ (respectively, $G^O(k)[C]$ and $G^O(k)[V \setminus C]$) provide tight examples of cobipartite graphs whose LD- and LTD-numbers (respectively, OD- and OTD-numbers) attain the logarithmic lower bound in Corollary~\ref{cor:log-LB}. Moreover, for $k \ge 2$, the graphs $G^L(k)$ (respectively, $G^O(k)$) with a complete edge set on the subgraph $G^L(k)[C]$ (respectively, $G^O(k)[C]$) and an empty edge set on $G^L(k)[V \setminus C]$ (respectively, $G^O(k)[V \setminus C]$) provide tight examples of split graphs whose LD- and LTD-numbers (respectively, OD- and OTD-numbers) attain the logarithmic lower bound in Corollary~\ref{cor:log-LB}. Finally, for $k \ge 2$, the graphs $G^I(k)$ with an empty edge set on the subgraph $G^L(k)[C]$ and a complete edge set on $G^L(k)[V \setminus C]$ provide tight examples of split graphs whose ID-numbers attain the logarithmic lower bound in Corollary~\ref{cor:log-LB}.

\begin{table}[ht]
\begin{center}
\begin{tabular}{ r || c | c || c | c || c}
 Graph family & LD & LTD & ID & OD & OTD \\ \hline
 Bipartite graphs & TIGHT & OPEN  & TIGHT & OPEN & OPEN \\
 Cobipartite graphs & TIGHT & TIGHT & OPEN & TIGHT & TIGHT \\
 Split graphs & TIGHT & TIGHT & TIGHT & TIGHT & TIGHT \\
\end{tabular}
\end{center}
\caption{
Immediate graph families some of whose members are tight examples of graphs with X-numbers attaining the logarithmic lower bounds in Corollary~\ref{cor:summary} for $X \in \{LD, LTD, ID, OD, OTD\}$.}
\label{tab_log LB_tight}
\end{table}

Apart from the tight examples given in Table~\ref{tab_log LB_tight}, our future work includes  the analysis of graph classes for which the logarithmic lower bounds hold as tight examples, or for which graph classes even stronger bounds can be established. Moreover, it would be interesting to provide similar constructions of extremal graphs for other or more general separation-domination problems.


\small
\begin{thebibliography}{11}

\bibitem{BCHL_2005}
Bertrand, N., Charon, I., Hudry, O., Lobstein, A.
\newblock \emph{1-identifying codes on trees.} 
\newblock The Australasian Journal of Combinatorics
\textbf{31}:~21--35 (2005). 

\bibitem{CW_ISCO2024} 
Chakraborty, D., Wagler, A.K.
\newblock \emph{Open-Separating Dominating Codes in Graphs.}
\newblock In: Basu, A., Mahjoub, A.R., Salazar González, J.J. (eds) Combinatorial Optimization. ISCO 2024, Lecture Notes in Computer Science
\textbf{14594}:~137--151 (2024).

\bibitem{CW_ODarxiv2024} 
Chakraborty, D., Wagler, A.K. 
\newblock \emph{Open-Separating Dominating Codes in Graphs.}
\newblock arXiv:2402.03015 [math.CO]
(2024).

\bibitem{CW_2024-FS} 
Chakraborty, D., Wagler, A.K.
\newblock \emph{On full-separating sets in graphs.}
\newblock arXiv:2407.10595 [math.CO]
(2024).

\bibitem{CHL_2003} Charon, I., Hudry, O., Lobstein, A. 
\newblock \emph{Minimizing the size of an identifying or locating-dominating code in a graph is {NP}-hard.}
\newblock Theoretical Computer Science
\textbf{290}:~2109--2120 (2003).


\bibitem{CSS_1987} 
Colburn, C., Slater, P.J., Stewart, L.K.
\newblock \emph{Locating-dominating sets in series-parallel networks.}
\newblock Congressus Numerantium
\textbf{56}:~135--162 (1987).


\bibitem{HHH_2006} Haynes, T.W., Henning, M. A., Howard, J.
\newblock \emph{Locating and total-dominating sets in trees.}
\newblock Discrete Applied Mathematics
\textbf{154}:~1293--1300 (2006).

\bibitem{HR_2012} Henning, M.A., Rad, J.R.
\newblock \emph{Locating-total domination in graphs.}
\newblock Discrete Applied Mathematics
\textbf{160}:~1986--1993 (2012).

\bibitem{HLR_2002} Honkala, I., Laihonen, T., Ranto, S.
\newblock \emph{On strongly identifying codes.}
\newblock Discrete Mathematics
\textbf{254}:~191--205 (2002).

\bibitem{KCL_1998}
M. G. Karpovsky, K. Chakrabarty, L. B. Levitin.
\newblock \emph{On a new class of codes for identifying vertices in graphs.}
\newblock IEEE Transactions on Information Theory \textbf{44}:~599--611 (1998).

\bibitem{Lobstein_Lib}
Jean, J., A. Lobstein.
\newblock \emph{Watching systems, identifying, locating-dominating and dis\-cri\-mi\-na\-ting codes in graphs.}\\[-8mm] 
%%%http://www.infres.enst.fr/lobstein/debutBIBidetlocdom.pdf
\newblock
\begin{verbatim}
https://www.lri.enst.fr/~lobstein/debutBIBidetlocdom.pdf
\end{verbatim}
%\bibitem{JL_lib} Jean, J., Lobstein, A.: Watching systems, identifying, locating-dominating and discriminating codes in graphs, https://dragazo.github.io/bibdom/main.pdf

\bibitem{M_2006} Moncel, J.
\newblock \emph{On graphs on $n$ vertices having an identifying code of cardinality $\log2(n + 1)$.}
\newblock Discrete Applied Mathematics
\textbf{154} (14):~2032 --2039 (2006).

\bibitem{S_1988}
P. J. Slater. 
\newblock \emph{Dominating and reference sets in a graph.}
\newblock Journal of Mathematical and Physical Sciences
\textbf{22}:~445--455 (1988).

\bibitem{SS_2010}
S. J. Seo, P. J. Slater.
\newblock \emph{Open neighborhood locating-dominating sets.}
\newblock Australasian Journal of Combinatorics \textbf{46}:~109--119 (2010).

\bibitem{UTS_2004}
Ungrangsi, R., Trachtenberg, A. and Starobinski, D.
\newblock \emph{An Implementation of Indoor Location Detection Systems Based on Identifying Codes.}
\newblock "Intelligence in Communication Systems", Springer Berlin Heidelberg,
175--189 (2004).


\end{thebibliography}
\end{document}